\documentclass[11pt, a4paper]{article}
\usepackage{latexsym}
\usepackage{indentfirst}
\usepackage{graphicx}
\usepackage{amsmath}
\usepackage{amssymb}
\begin{document}

\title{On the Computation of the Third Order Terms of the Series Defining the
Center Manifold for Systems of Delay Differential Equation}
\author{Anca - Veronica Ion \thanks{"Gh. Mihoc-C. Iacob" Institute of Mathematical Statistics and
Applied Mathematics of the Romanian Academy,
13, Calea 13 Septembrie, Bucharest, Romania, \newline anca-veronica.ion@ima.ro}}
\date{}
\maketitle

\begin{abstract} When studying a general system of delay differential equation with a single
constant delay, we encounter a certain lack of uniqueness in determining the coefficient of one of
the third order terms of the series defining the center manifold. We solve this problem by
considering a perturbation of the considered problem, perturbation that allows us to remove the
singularity.
The result generalizes a similar result obtained for scalar differential equations (J. Dyn.
Diff. Eqns., 24/2012).
\end{abstract}

\section{Introduction}

We consider an equation of the form
\begin{equation}\label{eq}\dot{x}(t)=f(x(t),x(t-r)),
\end{equation}
where $x:[-r,T)\rightarrow \mathbb{R}^n$, $f:D\subset \mathbb{R}^{2n}\rightarrow \mathbb{R}^n$, $(0,0)\in D\subset\mathbb{R}^{2n},$ $D$
open, $f(0,0)=0,$ $r>0.$ We assume $f\in C^{k}(D),\,k > 3$. We attach to this equation an initial condition of the form
\begin{equation}\label{ci} x(s)=\phi(s),\,s\in [-r,0],
\end{equation}
with $\phi: [-r,0]\rightarrow \mathbb{R}^n$ a continuous function.

Since $x(t)=0$ is a solution (an equilibrium) of \eqref{eq}, we consider the linearized of the
equation around this solution
\begin{equation}\label{lin-eq}\dot{x}(t)=\mathrm{A}x(t)+\mathrm{B}x(t-r),
\end{equation}
where $\mathrm{A},\,\mathrm{B}\in \mathcal{M}_{n\times n}(\mathbb{R}).$ Thus, eq. \eqref{eq} may
be written as
\begin{equation}\label{s-eq}\dot{x}(t)=\mathrm{A}x(t)+\mathrm{B}x(t-r)+\widehat{f}(x(t),\,x(t-r)),
\end{equation}
 $\widehat{f}$ being the nonlinear part of $f$.

The characteristic equation attached to \eqref{lin-eq} is
\begin{equation}\label{ch-eq}\det(\lambda I-\mathrm{A}-e^{-\lambda r}\mathrm{B})=0.
\end{equation}
\textbf{Hypothesis H.} \textit{We assume that \eqref{ch-eq} has a pair of pure imaginary
complex conjugate solutions $\lambda_{1,2}=\pm \omega i,$ with
$\omega>0,$ and all other eigenvalues have negative real part. We also assume that the eigenvalues
$\lambda_{1,2}$ are simple, (i.e. each of them has an one-dimensional corresponding eigenspace).}

\vspace{0.3cm}

Under this hypothesis, it is known that an invariant manifold, called \emph{center manifold},
exists in a certain function space that will be defined below. The dimension of the invariant
manifold is two, it is the graph of a function depending on two complex conjugated scalar variables.
The reduction of the problem to this manifold
allows some bifurcations studies to be performed, when the problem depends also upon some
parameters. 

In order to find an approximation of the reduced to the center manifold problem, an approximation
of the center manifold is required. This is obtained by considering the series of powers of the
function defining the center manifold, and by computing  (using a method that will be recalled in
Section 2), the first coefficients of this series.

When one arrives to the third order terms, one sees that a singularity occurs, i.e.  the system of
algebraical equations determining a certain value of the coefficient of $z^2\overline{z}$ has
singular matrix. In this situation either the problem has an infinity of solutions, or it has no
solution at all.

We show that our problem has an infinity of solutions (Section 3) and we present a method to
select a significant solution from among them (Section 4). A similar problem was solved for the
case of a single differential delay equation ($n=1$) in \cite{AVI-12}.

We must add that the knowledge of the third order terms of the center manifold series is necessary
for the study of degenerate Hopf bifurcation (when the problem presents one varying parameter) or
Bautin bifurcation (when the problem presents two varying parameters) - see, for a scalar
equation, \cite{AVI-13}.

\section{Preliminaries}

In this section the center manifold is defined and the method of approximating it is recalled. These are basic well-known
ideas, but we recall them for the completeness of the paper.

\subsection{The center manifold}

This subsection relies on  \cite{Ha}, \cite{HaL}, \cite{Far}.  In order to define the center manifold, we need the space
$$\mathcal{B}=\left\{\psi:[-r,0]\mapsto \mathbb{R}^n\,\,|\,\,\, \psi\,
\mathrm{\, continuous\, on\,}[-r,0]\, \right\},$$ and its
complexification $\mathcal{B}_C=\mathcal{B}+i\mathcal{B}$, both with sup norm.

The eigenfunctions corresponding to $\lambda_{1,2}=\pm \omega i$ (see hypothesis  \textbf{H}) are
the elements
of $\mathcal{B}_C$ given by
$\varphi_{1,2}(s)=\varphi_{1,2}(0)e^{\pm\omega is},$ where
$\varphi_{1,2}(0)$ are solutions of
\begin{equation}\label{ev}
(\pm\omega iI-\mathrm{A} -e^{\mp\omega i r}\mathrm{B})\varphi_{1,2}(0)=0.\end{equation}
We consider the subspace of $\mathcal{B}_C$ defined by
$\mathcal{M}=Span\{\varphi_{1},\,\varphi_{2}\}$ ($\varphi_{1,2}$ are column vectors). 
In \cite{Ha}, \cite{HaL} a projector $\mathcal{P}:\mathcal{B}_C\rightarrow \mathcal{M}$
 is defined, and by setting \break $\mathcal{N}=(I-\mathcal{P})\mathcal{B}_C,$ we have
 $\mathcal{B}_C=\mathcal{M}\oplus
\mathcal{N}$. We review here the construction of this projector.
The adjoint equation (associated to the linear equation
\eqref{lin-eq}) is defined as
\begin{equation}\dot{y}(s)=-\mathrm{A}y(s)-\mathrm{B}y(s+r).
\end{equation}
The corresponding characteristic equation is
\[\det(\lambda \mathrm{I}+\mathrm{A}+\mathrm{B}e^{\lambda r})=0,
\]
that, obviously has, together with eq. \eqref{ch-eq},
the solutions $\pm \omega i.$\\
 The corresponding eigenfunctions
 are $\psi_{1}(\zeta)=\psi_{1}(0)e^{-\omega i\zeta},\;\psi_{2}(\zeta)=\psi_{2}(0)e^{\omega
 i\zeta},$\\$\zeta\in[0,r]$ where $\psi_j(0)$ are row vectors, solutions of the vectorial
 equations
\begin{equation}\label{ev-adj}\psi_{1,2}(0)(\mp\omega i
\mathrm{I}+\mathrm{A}+\mathrm{B}e^{\mp\omega i r} )=0.
\end{equation}

The following bilinear form is defined (\cite{Ha},
\cite{HaL}) on
$C([0,r],\mathbb{C}^n)\times \mathcal{B}_C$
\begin{equation}\label{bf}\langle\psi,\varphi\rangle=\psi(0)\varphi(0)+
\int_{-r}^0\psi(\zeta+r) \mathrm{B}\varphi(\zeta)d\zeta,
\end{equation}
where $\psi$ is a row vector and $\varphi$ is a
column vector.

Linear combinations of the functions $\psi_j, \;j=1,2,$ denoted
by \break $\Psi_i,\,i=1,2$ are then constructed, such that $\langle \Psi_i,
\varphi_j \rangle=\delta_{ij}.$ For this, the matrix $E=(e_{ij})_{1\leq i,j\leq 2}$, with $e_{ij}=\langle \psi_i, \,\varphi_j
\rangle$, is computed (see the Appendix). We obtain
$\Psi_1(\zeta)=\Psi_1(0)e^{-\omega i \zeta}$, with
$$\Psi_1(0)=\frac{1}{\psi_1(0)\varphi_1(0)+\psi_1(0)
 \mathrm{B}\varphi_1(0)e^{-\omega
ri}r}\psi_1(0).$$

The projector defined (\cite{Ha}, \cite{HaL}) on $\mathcal{B}_{C}$,
 with values in $\mathcal{M}$  is given by
$$\mathcal{P}(\phi)=\langle\Psi_1,\phi\rangle\varphi_1
+\langle\Psi_2,\phi\rangle\varphi_2,\,\,\phi\in\mathcal{B}_C.$$

\vspace{0.2cm}

In the hypothesis \textbf{H} above, it is proved \cite{HaL} that  a
local invariant manifold, called the \textit{center manifold}, exists, that
is a smooth manifold, tangent to the space $\mathcal{M}$ at the
point $x=0$  and being the graph of a function
$w(\cdot)$ defined on a neighborhood of zero in
$\mathcal{M}$ with values in $\mathcal{N}$. Since all eigenvalues different of $\lambda_{1,2}$ have strictly negative real part, the trajectory of any point that does not lie on the center manifold approaches, when $t \rightarrow \infty,$ the center manifold.

In order to approximate the center manifold, the function \break $w(u,\,\overline{u}):=w(u\,\varphi_1+\overline{u}\,\varphi_2)$ is written as a series of powers $u$ and $\overline{u}$,
\begin{equation}\label{cm}  w(u,\overline{u})=\sum_{i+j\geq 2}\frac{1}{i!j!}w_{i,j}u^i\overline{u}^j,\end{equation}
 where for each $i,j,$ $w_{i,j}\in \mathcal{B}_C$, hence it is a vectorial function on $[-r,0]$.

\subsection{Computation of the coefficients $w_{i,j}$}

The vectorial functions $w_{i,j}$ of \eqref{cm} are determined by using the
following relation  \cite{WWPOG}, \cite{MNO}, \cite{AVI} that is a consequence of the invariance of the
center manifold.

\begin{equation}\label{gendiffeq}\frac{\partial}{\partial
s}\sum_{j+k\geq2}\frac{1}{j!k!}w_{j,k}(s)u^{j}\overline{u}^{k}=
\sum_{j+k\geq2}\frac{1}{j!k!}g_{j,k}u^{j}\overline{u}^{k}\varphi_1(s)+
\end{equation}
\[+
\sum_{j+k\geq2}\frac{1}{j!k!}\overline{g}_{j,k}\overline{u}^{j}u^{k}\varphi_2(s)+\frac{\partial}{\partial
t}\sum_{j+k\geq2}\frac{1}{j!k!}w_{j,k}(s)u^{j}\overline{u}^{k}.
\]
From this, differential equations for each $w_{j,k}$ are obtained, by matching the terms of degree
$j$ in $u$ and of degree $k$ in $\overline{u}.$
To determine the integration constants, the following relation is used:

\begin{equation}\label{conddiffeq}
\frac{d}{dt}\sum_{j+k\geq2}\frac{1}{j!k!}w_{j,k}(0)u^{j}\overline{u}^{k}+
\sum_{j+k\geq2}\frac{1}{j!k!}g_{j,k}u^{j}\overline{u}^{k}\varphi_1(0)+
\sum_{j+k\geq2}\frac{1}{j!k!}\overline{g}_{j,k}\overline{u}^{j}u^{k}\varphi_2(0)=
\end{equation}
\[=\mathrm{A}\sum_{j+k\geq2}\frac{1}{j!k!}w_{j,k}(0)u^{j}\overline{u}^{k}+\mathrm{B}\sum_{j+k\geq2}\frac{1}{j!k!}w_{j,k}(-r)u^{j}\overline{u}^{k}+
\sum_{j+k\geq2}\frac{1}{j!k!}f_{j,k}u^{j}\overline{u}^{k}.
\]

\section{The system for  $w_{2,1}(-r),\,\,w_{2,1}(0)$}

By matching, in \eqref{gendiffeq} and \eqref{conddiffeq}, the terms that contain
$u^2\overline{u}$,  we get the system of
differential  equations and the set of conditions for $w_{2,1}.$ These are \newpage

\[\frac{d }{ds}w_{2,1}(s)=\omega i w_{2,1}(s)+g_{2,1}\varphi_{1}(0)e^{\omega i
s}+\overline{g}_{1,2}\varphi_{2}(0)e^{-\omega is}+2w_{2,0}(s)g_{1,1}+\]
\[\quad\quad\quad\quad+w_{1,1}(s)g_{2,0}+2w_{1,1}(s)\overline{g}_{1,1}
+w_{0,2}(s)\overline{g}_{0,2},
\]
(that is a system of $n$ equations) and the $n$ conditions:
\[\omega iw_{2,1}(0)
+2w_{2,0}(0)g_{1,1}+w_{1,1}(0)g_{2,0}+2w_{1,1}(0)\overline{g}_{1,1}
+w_{0,2}(0)\overline{g}_{0,2}+g_{2,1}\varphi_{1}(0)+\overline{g}_{1,2}\varphi_{2}(0)=\]
\[
=\mathrm{A}w_{2,1}(0)+ \mathrm{B} w_{2,1}(-r)+f_{2,1}.
\]

From these  we obtain the system of equations for $w_{21}(0)$  and
$w_{21}(-r)$ (vectorial equations):

\begin{equation}\label{alg-sys1}-e^{-\omega ir}w_{2,1}(0)+w_{2,1}(-r)=-g_{2,1}re^{-\omega
ir}\varphi_{1}(0)+\frac{i}{2\omega} \overline{g}_{1,2}(e^{\omega ir}-e^{-\omega
ir})\varphi_{2}(0)-\end{equation} \[-2g_{1,1}e^{-\omega
ir}\int_{-r}^0w_{2,0}(\theta)e^{-\omega i
\theta}d\theta-(g_{2,0}+2\overline{g}_{1,1})e^{-\omega
ir}\int_{-r}^0 w_{1,1}(\theta)e^{-\omega i\theta}d\theta -\]
\[-\overline{g}_{0,2}e^{-\omega ir}\int_{-r}^0
w_{0,2}(\theta)e^{-\omega i\theta}d\theta,
\]
\medskip
\[-(\omega i\mathrm{I}-\mathrm{A})w_{2,1}(0)+
\mathrm{B}w_{2,1}(-r)=g_{2,1}\varphi_{1}(0)+\overline{g}_{1,2}\varphi_{2}(0)-f_{2,1}+2g_{1,1}w_{2,0}(0)+
\]
\begin{equation}\label{alg-sys2}\quad\quad\quad\quad+
(g_{2,0}+2\overline{g}_{1,1})w_{1,1}(0)
+\overline{g}_{0,2}w_{0,2}(0).\end{equation}

The matrix of this system
is the $2n\times 2n$ matrix (below, `` $\mathrm{I}$ ''  is the $n\times n$ identity matrix)
\[
    \left(
\begin{array}{cc}
 -e^{-\omega ir}\mathrm{I}& \mathrm{I}\\
       -(\omega i\mathrm{I-A})   & \mathrm{B} \end{array}
     \right).
\]
 We denote the right hand side of \eqref{alg-sys1} by
$R_1$ and the right hand side of \eqref{alg-sys2} by $R_2.$ We multiply the vectorial equation
\eqref{alg-sys1} to the left  by $\mathrm{B}$ and, by subtracting the second equation from the
first we get the system
\begin{equation}\label{new-sys1}(\omega i \mathrm{I-A-B}e^{-\omega
ir})w_{21}(0)=\mathrm{B}R_1-R_2.
\end{equation}
Since $\omega i$ is an eigenvalue for our problem, the determinant of the matrix of the above
system is equal to zero. Hence, at this moment, we don't know whether the system has solutions or
not. The following result solves this problem.

\medskip

\textbf{Proposition 3.1.} \textit{The system of equations \eqref{new-sys1} has solutions.}

\smallskip

\textbf{Proof.} Since $-\omega i$ is a simple eigenvalue, the matrix $$\mathrm{M}=\omega i
\mathrm{I}-\mathrm{A}-\mathrm{B}e^{-\omega ir}$$ has an one-dimensional kernel. The same is true
for the matrix \\$(\omega i \mathrm{I}-\mathrm{A}-\mathrm{B}e^{-\omega ir})^{T}$.
The kernel of this latter matrix is spanned by $\Psi_1(0)^{T},$ since we have
\[\Psi_1(0)(\omega i \mathrm{I}-\mathrm{A}-\mathrm{B}e^{-\omega ir})=0.
\]

For a $n\times n$ matrix $M$ the equality
$\mathcal{R}(\mathrm{M})=(\mathcal{N}(\mathrm{M}^{T}))^{\bot}$ holds (where $\mathcal{R}$ is the
range and $\mathcal{N}$ is the null space of the matrices in brackets). Hence, in order to prove
that system \eqref{new-sys1} has solutions, we show that $\mathrm{B}R_1-R_2$ is orthogonal on
$\mathcal{N}\left((\omega i \mathrm{I}-\mathrm{A}-\mathrm{B}e^{-\omega ir})^{T}\right), $ that is
\begin{equation}\label{BR}\Psi_1(0)(\mathrm{B}R_1-R_2 )=0.\end{equation}

We
prove this by showing that the following relations hold:
$$\Psi_1(0)\left(-g_{2,1}re^{-\omega ir}\mathrm{B}\varphi_1(0)-g_{2,1}\varphi_{1}(0)+f_{2,1}
\right)=0; \eqno{(R1a)} $$
$$\Psi_1(0)\left(\frac{i}{2\omega}
\overline{g}_{1,2}\mathrm{B}\overline{\varphi}_1(0)(e^{\omega ir}-e^{-\omega
ir})-\overline{g}_{1,2}\overline{\varphi}_{1}(0) \right)=0; \eqno{(R1b)}$$
$$\Psi_1(0)\left(-2g_{1,1}e^{-\omega
ir}\mathrm{B}\int_{-r}^0w_{2,0}(\theta)e^{-\omega
i\theta}d\theta-2g_{1,1}w_{2,0}(0)\right)=0; \eqno{(R2)}$$
$$\Psi_1(0)\left(-(g_{2,0}+2\overline{g}_{1,1})e^{-\omega ir}\mathrm{B}\int_{-r}^0
w_{1,1}(\theta)e^{-\omega
i\theta}d\theta-(g_{2,0}+2\overline{g}_{1,1})w_{1,1}(0)\right)=0;
\eqno{(R3)}$$
$$\Psi_1(0)\left(\overline{g}_{0,2}e^{-\omega ir}-\mathrm{B}\int_{-r}^0
w_{0,2}(\theta)e^{-\omega
i\theta}d\theta-\overline{g}_{0,2}w_{0,2}(0)\right)=0. \eqno{(R4)}$$
\vspace{0.5cm}

We prove these relations one by one.\\
\textbf{(\textit{R}1\textit{a}).} This can be written as

\[-g_{2,1}re^{-\omega ir}\Psi_1(0) B\varphi_1(0)-g_{2,1}\Psi_1(0)\varphi_{1}(0)+g_{2,1}=0.
\]
If $g_{2,1}=0,$ then the relation is satisfied, if $g_{2,1}\neq 0,$ we divide by $g_{2,1}$ and
obtain
\[-re^{-\omega ir}\Psi_1(0)\mathrm{B}\varphi_1(0)-\Psi_1(0)\varphi_{1}(0)+1=0,
\]
that we rewrite as
\[\Psi_1(0)\varphi_{1}(0)+\Psi_1(0)\mathrm{B}\varphi_1(0)
\int_{-r}^{0}e^{-\omega i(\zeta +r)}e^{\omega i \zeta}d\zeta=1.
\]
But this is the true relation
\[\langle \Psi_1, \varphi_1\rangle =1.
\]
Relation \textbf{(R1\textit{b})} leads, after dividing it by $\overline{g}_{1,2}$ (assumed
non-null)  to
$$\Psi_1(0)\overline{\varphi}_{1}(0)+\frac{i}{2\omega}
\Psi_1(0)\mathrm{B}\overline{\varphi}_1(0)(e^{-\omega ir}-e^{\omega
ir})=0 $$
that can be written as
$$\Psi_1(0)\overline{\varphi}_{1}(0)+
\Psi_1(0)\mathrm{B}\overline{\varphi}_1(0)\int_{-r}^0(e^{-\omega i(\zeta+r)} e^{-\omega i
\zeta})d\zeta=0 $$
this last being equivalent to
$$\langle \Psi_1,\,\overline{\varphi}_1\rangle=0,
$$
that holds true.

\textbf{Relations (\textit{R}2)-(\textit{R}4).} If $g_{1,1} =0,$
then relation (\textit{R}2) holds. We assume $g_{1,1} \neq 0.$
We also assume
$g_{2,0}+2\overline{g}_{1,1}\neq 0$ and $\overline{g}_{0,2}\neq 0$.

After dividing with the assumed non-zero coefficients,
 relations (\textit{R}2)-(\textit{R}4) have the form:
$$\Psi_1(0)\left(-e^{-\omega
ir}\mathrm{B}\int_{-r}^0w_{jk}(\theta)e^{-\omega
i\theta}d\theta-w_{jk}(0)\right)=0; $$
that can be written as
$$\Psi_1(0) w_{jk}(0)+\int_{-r}^0\Psi_1(0)e^{-\omega
i(\theta+r)} \mathrm{B}w_{jk}(\theta)d\theta=0
$$
($j,k\geq 0,\,j+k=2$) that  is
\[\langle \Psi_1, w_{j,k}\rangle =0,
\]
that holds true, since each of the functions $w_{j,k}$ belongs to the
complementary of $\mathcal{M}.$ Hence all relations (\textit{R}2) -
(\textit{R}4) are proved.

\smallskip

By adding the five relations proved we obtain \eqref{BR} and the
conclusion of our Proposition follows. $\Box$

\section{How to compute $w_{2,1}(0)$}

Now, in order to solve the system for $w_{2,1}(0)$, we use the following ideas. We take the vector
$v_1=\Psi_1(0)^T/|\Psi_1(0)|,$ and we complete it to a  orthonormal basis in $\mathbb{C}^n$,
$\mathcal{B}=\{v_1,\,v_2,\,...,v_n\}$ (column vectors).

We write the system $\mathrm{M}w_{2,1}(0)=\mathrm{B}R_1-R_2$ in the basis $\mathcal{B}$.

Thus,we set $x_j=w_{2,1}(0)\cdot v_j$ (the coordinates of $w_{2,1}(0)$ with respect to
$\mathcal{B}$), and we obtain the system in the new basis (the ``$\cdot$'' represents the scalar
product in the space of $n$ dimensional column vectors):

$$\left\{\begin{array}{lll}
(v_1\cdot \mathrm{M}v_1)x_1+(v_1\cdot \mathrm{M}v_2)x_2+...+(v_1\cdot
\mathrm{M}v_n)x_n&=&v_1\cdot(\mathrm{B}R_1-R_2),\\
(v_2\cdot \mathrm{M}v_1)x_1+(v_2\cdot \mathrm{M}v_2)x_2+...+(v_2\cdot
\mathrm{M}v_n)x_n&=&v_2\cdot(\mathrm{B}R_1-R_2),\\
                                                      &\vdots &                     \\
(v_n\cdot \mathrm{M}v_1)x_1+(v_n\cdot \mathrm{M}v_2)x_2+...+(v_n\cdot
\mathrm{M})v_nx_n&=&v_n\cdot(\mathrm{B}R_1-R_2).\\
\end{array}\right.$$

Since $\Psi_1(0)\mathrm{M}=0$ (here $0$ is a row vector), we have $v_1\cdot \mathrm{M}v_j=0$ for
any $j=1,...,n$ and since $v_1\cdot(\mathrm{B}R_1-R_2)=0,$ the first equation is identically
satisfied, hence we can not use it.

We assumed that the eigenvalues $\pm \omega i$ are simple, hence the rank of the matrix
$\mathrm{M}$ is $n-1$, and the system of $n-1$ equations obtained by removing the first equation
can be solved with respect to $n-1$ variables. Thus only the first equation is responsible for the
singularity of the problem.

\vspace{0.2cm}

In order to remove this singularity, as in the scalar case, we  consider a perturbation of our
problem, of the form

\begin{equation}\label{pert-eq}\dot{x}(t)=\mathrm{A}_\epsilon x(t)+\mathrm{B}_\epsilon
x(t-r)+\widehat{f}(x(t),\,x(t-r)),
\end{equation}
where $\mathrm{A}_\epsilon, \,\mathrm{B}_\epsilon$ depend smoothly enough on
$\epsilon>0,$ $\displaystyle\lim_{\epsilon\searrow 0}\mathrm{A}_\epsilon =
\mathrm{A},\, \lim_{\epsilon\searrow 0}\mathrm{B}_\epsilon = \mathrm{B},$ and they are chosen
such
that, for small enough $\epsilon$,  the linearized problem attached
to \eqref{pert-eq} admits the eigenvalues $\lambda_{\epsilon
1,2}=\mu_\epsilon\pm \omega_\epsilon i,$ with $\mu_\epsilon>0,$
while all other eigenvalues have negative real part.

It follows that $\displaystyle\lim_{\epsilon\searrow 0}\mu_\epsilon = 0,\,
\lim_{\epsilon\searrow 0}\omega_\epsilon = \omega$.

Since $\mu_\epsilon>0,$ the problem \eqref{pert-eq}
admits an unstable manifold, tangent to the space
$\mathcal{M}_\epsilon$ spanned by the two eigenfunctions
corresponding to the two eigenvalues $\lambda_{\epsilon 1,2},$ i.e.
$\varphi_{\epsilon 1,2}(s)=\varphi_{\epsilon 1,2}(0)e^{(\mu_\epsilon\pm i\omega_\epsilon)s}.$

This unstable manifold is the graph of a function $w_{\epsilon}$
defined on $\mathcal{M}_\epsilon$ and with values in a subspace
$\mathcal{N}_\epsilon$ of $\mathcal{B}_C,$ such that
$\mathcal{B}_C=\mathcal{M}_\epsilon\oplus \mathcal{N}_\epsilon$.

Now, as in Section 2 for the
non-perturbed problem, the following mathematical objects are constructed \cite{Far}:

-the adjoint of the linear equation and its eigenfunctions  i.e. the
functions $\psi_{\epsilon 1}(s)=\psi_{\epsilon 1}(0)e^{-\lambda_\epsilon
s},\,\psi_{\epsilon 2}(s)=\overline{\psi}_{\epsilon 1}(0)e^{-\overline{\lambda}_\epsilon
s},\,s\in
[0,r]$;

-the corresponding bilinear form, denoted  also by $\langle\,
\cdot\, ,\,\cdot\,\rangle$ given by
\[\langle \psi,\varphi \rangle=\psi(0)\varphi(0)+
\int_{-r}^0\psi(\zeta+r)\mathrm{B}_{\epsilon}\varphi(\zeta)d\zeta.
\]

We find the functions $\Psi_{\epsilon 1}, \,\Psi_{\epsilon 2},$ such that $\langle \Psi_{\epsilon
i}\,,\, \varphi_{\epsilon j} \rangle =\delta_{i,j}$ (see the Appendix).  The projector
$\mathcal{P}_\epsilon:\mathcal{B}_C\rightarrow
\mathcal{M}_\epsilon,$  is defined by
$$\mathcal{P}_\epsilon(\phi)=\langle \Psi_{\epsilon 1},\phi\rangle \varphi_1+\langle
\Psi_{\epsilon 2},\phi\rangle \varphi_2.$$ We set
$\mathcal{N}_\epsilon=(I-\mathcal{P}_\epsilon)(\mathcal{B}_C)$. Note that
\[\langle \Psi_{\epsilon j}, w_\epsilon \rangle=0,\,\,j=1,2.
\]
 since $w_\epsilon$ takes values in $\mathcal{N}_\epsilon$.

The problem reduced to the unstable manifold is \cite{Far}
\begin{equation}\frac{dv}{dt}=\lambda_{\epsilon}v+\Psi_{\epsilon
1}(0)\widetilde{f}(v\varphi_{\epsilon
1}+\overline{v}\varphi_{\epsilon 2}+w_\epsilon(v\varphi_{\epsilon
1}+\overline{v}\varphi_{\epsilon 2})).
\end{equation}

We consider the function
$$w_\epsilon(v,
\overline{v}):=w_\epsilon(v\varphi_{\epsilon
1}+\overline{v}\varphi_{\epsilon 2})$$ and  we write:
\[  w_\epsilon(v,\overline{v})=\sum_{i+j\geq 2}\frac{1}{i!j!}w_{\epsilon i,j}v^i\overline{v}^j.
\]
As in the non-perturbed case, the coefficients  $w_{\epsilon i,j}$ are found by solving
differential equations coming from  relations similar to
\eqref{gendiffeq}, \eqref{conddiffeq} with $w_{j,k}, \,g_{j,k},...$
replaced by  $w_{\epsilon j,k}, \,g_{\epsilon j,k},...$.

\vspace{0.2cm}

The equation for $w_{\epsilon 2,1}$ is

\[\frac{d w_{\epsilon2,1}(s)}{ds}=(2\lambda_\epsilon+\overline{\lambda}_\epsilon) w_{\epsilon
2,1}(s)+g_{\epsilon 2,1}\varphi_{\epsilon1}(0)e^{\lambda_\epsilon
s}+\overline{g}_{\epsilon 1,2}\overline{\varphi}_{\epsilon1}(0)e^{\overline{\lambda}_\epsilon s}
+2w_{\epsilon 2,0}(s)g_{\epsilon 1,1}+\]
\[+w_{\epsilon 1,1}(s)g_{\epsilon 2,0}+2w_{\epsilon 1,1}(s)\overline{g}_{\epsilon 1,1}
+w_{\epsilon 0,2}(s)\overline{g}_{\epsilon 0,2},
\]
while the condition to determine the integration constant is
\[(2\lambda_\epsilon+\overline{\lambda}_\epsilon) w_{\epsilon 2,1}(0)
+2w_{\epsilon 2,0}(0)g_{11}+w_{\epsilon
1,1}(0)g_{\epsilon2,0}+2w_{\epsilon1,1}(0)\overline{g}_{\epsilon1,1}
+w_{\epsilon0,2}(0)\overline{g}_{\epsilon0,2}+\]
\[+g_{\epsilon2,1}\varphi_{\epsilon1}(0)+
\overline{g}_{\epsilon1,2}\overline{\varphi}_{\epsilon1}(0)
=\mathrm{A}_{\epsilon}w_{\epsilon2,1}(0)+
\mathrm{B}_{\epsilon}w_{\epsilon2,1}(-r)+f_{\epsilon2,1}.
\]
By integrating the differential equation above between $-r$ and $0$,
we find the system of equations for
$w_{\epsilon2,1}(-r),\,w_{\epsilon2,1}(0):$

\begin{equation}\label{sist-eps1}-e^{-(2\lambda_\epsilon+\overline{\lambda}_\epsilon)r}w_{\epsilon2,1}(0)+w_{\epsilon2,1}(-r)=
\frac{-1}{\lambda_\epsilon+\overline{\lambda}_\epsilon}g_{\epsilon2,1}\varphi_{\epsilon1}(0)(
e^{-\lambda_\epsilon
r}-e^{-(2\lambda_\epsilon+\overline{\lambda}_\epsilon)r})-\;\;\;
\end{equation}
\[-\frac{1}{2\lambda_\epsilon}\overline{g}_{\epsilon1,2}\overline{\varphi}_{\epsilon1}(0)(e^{-\overline{\lambda}_\epsilon
r}-e^{-(2\lambda_\epsilon+\overline{\lambda}_\epsilon)r})-2g_{\epsilon1,1}e^{-(2\lambda_\epsilon+\overline{\lambda}_\epsilon)r}
\int_{-r}^0w_{\epsilon2,0}(\theta)e^{-(2\lambda_\epsilon+
\overline{\lambda}_\epsilon) \theta}d\theta-\]

\[-(g_{\epsilon2,0}+2\overline{g}_{\epsilon1,1})e^{-(2\lambda_\epsilon+\overline{\lambda}_\epsilon)r}
\int_{-r}^0
w_{\epsilon1,1}(\theta)e^{-(2\lambda_\epsilon+\overline{\lambda}_\epsilon)
\theta}d\theta -
\]
\[-\overline{g}_{\epsilon0,2}e^{-(2\lambda_\epsilon+\overline{\lambda}_\epsilon)r}\int_{-r}^0
w_{\epsilon0,2}(\theta)
e^{-(2\lambda_\varepsilon+\overline{\lambda}_\epsilon)
\theta}d\theta,
\]
\vspace{0.8cm}
\begin{equation}\label{sist-eps2}
\left(\mathrm{A}_{\epsilon}-2\lambda_\epsilon \mathrm{I}-\overline{\lambda}_\epsilon
\mathrm{I}\right)w_{\epsilon2,1}(0)+
\mathrm{B}_{\epsilon}w_{\epsilon2,1}(-r)=g_{\epsilon2,1}\varphi_1(0)+
\overline{g}_{\epsilon1,2}\overline{\varphi}_1(0)-f_{\epsilon2,1}+
\end{equation}
\[+2w_{\epsilon 2,0}(0)g_{11}+w_{\epsilon
1,1}(0)g_{\epsilon2,0}+2w_{\epsilon1,1}(0)\overline{g}_{\epsilon1,1}
+w_{\epsilon0,2}(0)\overline{g}_{\epsilon0,2}.
\]

\medskip
As in the scalar case, we can easily prove the following result.

\textbf{Proposition 4.1.} \emph{When $\epsilon\rightarrow 0,$ the
coefficients of system \eqref{sist-eps1}-\eqref{sist-eps2} tend to
the coefficients of system \eqref{alg-sys1}-\eqref{alg-sys2}.}
\smallskip

The proof repeats the ideas of the proof of the similar assertion from the scalar problem, in
\cite{AVI-12}, hence we skip it.\\

We now denote the right-hand side of the two equations by
$R_{\epsilon1},\,R_{\epsilon2}, $ respectively. As in the non-perturbed case, we multiply the
first vectorial equation at the left by $\mathrm{B}_{\epsilon}$ and we subtract the second
vectorial equation from the such modified first equation.

The matrix of the obtained $n\times n$ system is
\begin{equation}\label{Delta-eps}\mathrm{M}_\epsilon=-\mathrm{B}_{\epsilon}e^{-(2\lambda_\epsilon+\overline{\lambda}_\epsilon)r}-
\mathrm{A}_{\epsilon}+2\lambda_\epsilon\mathrm{I}+\overline{\lambda}_\epsilon\mathrm{I},\end{equation}
and the system has the form:
\begin{equation}\label{sysw0}\mathrm{M}_\epsilon w_{\epsilon2,1}(0)=\mathrm{B}_\epsilon
R_{\epsilon 1}-R_{\epsilon 2}.
\end{equation}

Its matrix, $\mathrm{M}_\epsilon,$ has determinant different of zero because otherwise the number
$2\lambda_\epsilon+\overline{\lambda}_\epsilon$ would be an
eigenvalue, with real part equal to $3\mu_\epsilon$, that
contradicts the fact that all eigenvalues have real part $\leq
\mu_\epsilon$.

Now we construct a basis, $\mathcal{B}_\epsilon$, of orthonormal vectors in $\mathbb{C}^n$,
containing the vector
$v_{\epsilon 1}=\Psi_{\epsilon 1}(0)^T/|\Psi_{\epsilon 1}(0)|$ (column vectors).

The vectors of this basis, $v_{\epsilon j},\,j=2,...,n$, should be chosen such that
$\lim_{\epsilon\rightarrow 0}v_{\epsilon j}=v_{j},$ where $v_{j}$ are the vectors of the basis
$\mathcal{B}$ from the non-perturbed case. We express the system of equations \eqref{sysw0} as
$$\left\{\begin{array}{lll}
(v_{\epsilon 1}\cdot \mathrm{M}_{\epsilon}v_{\epsilon 1})x_1+(v_{\epsilon 1}\cdot
\mathrm{M}_{\epsilon}v_{\epsilon 2})x_2+...+(v_{\epsilon 1}\cdot \mathrm{M}_{\epsilon}v_{\epsilon
n})x_n&=&v_{\epsilon 1}\cdot(\mathrm{B}_{\epsilon}R_{\epsilon 1}-R_{\epsilon 2}),\\
(v_{\epsilon 2}\cdot \mathrm{M}_{\epsilon}v_{\epsilon 1})x_1+(v_{\epsilon 2}\cdot
\mathrm{M}_{\epsilon}v_{\epsilon 2})x_2+...+(v_{\epsilon 2}\cdot \mathrm{M}_{\epsilon}v_{\epsilon
n})x_n&=&v_{\epsilon 2}\cdot(\mathrm{B}_{\epsilon}R_{\epsilon 1}-R_{\epsilon 2}),\\
                                                      &\vdots &                     \\
(v_{\epsilon n}\cdot \mathrm{M}_{\epsilon}v_{\epsilon 1})x_1+(v_{\epsilon n}\cdot
\mathrm{M}_{\epsilon}v_{\epsilon 2})x_2+...+(v_{\epsilon n}\cdot \mathrm{M}_{\epsilon}v_{\epsilon
n})x_n&=&v_{\epsilon n}\cdot(\mathrm{B}_{\epsilon}R_{\epsilon 1}-R_{\epsilon 2}).\\
\end{array}\right.$$

\vspace{0.2cm}

The way we chose the vectors ensures us that, when $\epsilon \rightarrow 0,$ the above system
tends to the one obtained for $\epsilon=0$.
Now we shall treat the first equation of this system. First we prove

\vspace{0.2cm}

\textbf{Proposition 4.2.} \textit{For any $\epsilon>0$, we have}
\[\Psi_{\epsilon 1}(0)\mathrm{M}_\epsilon=\mu_\epsilon
h_1(\epsilon),
\]
\[\Psi_{\epsilon 1}(0)(\mathrm{B}_{\epsilon}R_{\epsilon1}-R_{\epsilon2})=\mu_\epsilon
h_2(\epsilon),
\]
\textit{where the functions $h_j,\,j=1,2$ have finite limit for
$\epsilon\rightarrow 0.$  }

\vspace{0.2cm}

\textbf{Proof.} By using the fact that $\Psi_{\epsilon 1}(0)$ is an eigenvalue for the perturbed
problem, i.e.
\[\Psi_{\epsilon 1}(0)(\lambda_\epsilon\mathrm{I}-\mathrm{A}_\epsilon-\mathrm{B}_\epsilon
e^{-\lambda_\epsilon  r})=0,\]
we have
\[\Psi_{\epsilon 1}(0)\mathrm{M}_\epsilon=\Psi_{\epsilon
1}(0)(-\mathrm{B}_{\epsilon}e^{-(2\lambda_\epsilon+\overline{\lambda}_\epsilon)r}-
\mathrm{A}_{\epsilon}+2\lambda_\epsilon\mathrm{I}+\overline{\lambda}_\epsilon\mathrm{I})=\]
\[=\Psi_{\epsilon
1}(0)\left(-\mathrm{B}_{\epsilon}e^{-(2\lambda_\epsilon+\overline{\lambda}_\epsilon)r}+\mathrm{B}_\epsilon
e^{-\lambda_\epsilon r}+\lambda_\epsilon \mathrm{I}+\overline{\lambda}_\epsilon
\mathrm{I}\right)=
\]
\[=\Psi_{\epsilon 1}(0)\left(\mathrm{B}_\epsilon e^{-\lambda_\epsilon r}\left(1-e^{-2\mu_\epsilon
r}\right)+2\mu_\epsilon \mathrm{I}\right)=\]
\[=\Psi_{\epsilon 1}(0)\left(\mathrm{B}_\epsilon e^{-\lambda_\epsilon r}2\mu_\epsilon
r\left(1-\frac{2\mu_\epsilon r}{2!}+\frac{(2\mu_\epsilon
r)^2}{3!}+...\right)+2\mu_\epsilon \mathrm{I}\right)=
\]
\[=2\mu_\epsilon\Psi_{\epsilon 1}(0)\left[\mathrm{B}_\epsilon e^{-\lambda_\epsilon r}
r\left(1-\frac{2\mu_\epsilon r}{2!}+\frac{(2\mu_\epsilon r)^2}{3!}+...\right)+\mathrm{I}\right].
\]
By denoting
\[h_{1}(\epsilon)=2\Psi_{\epsilon 1}(0)\left[\mathrm{B}_\epsilon e^{-\lambda_\epsilon r}
r\left(1-\frac{2\mu_\epsilon r}{2!}+\frac{(2\mu_\epsilon r)^2}{3!}+...\right)+\mathrm{I}\right]
\]
and remarking that  $$\lim_{\epsilon\rightarrow 0}h_1(\epsilon)=2\Psi_{1}(0) \left(\mathrm{B}
e^{-\omega i r} r+\mathrm{I}\right),$$
we obtain the first assertion of the Proposition.

Now, in order to treat the term
$\mathrm{B}_{\epsilon}R_{\epsilon1}-R_{\epsilon2}$, inspired by the proof of
Proposition 3.1, we write
\[\Psi_{\epsilon
1}(0)(\mathrm{B}_{\epsilon}R_{\epsilon1}-R_{\epsilon2})=E_{\epsilon1a}+E_{\epsilon1b}+E_{\epsilon2}+E_{\epsilon3}+E_{\epsilon4},
\]
where
\[E_{\epsilon
1a}=\frac{-g_{\epsilon2,1}}{\lambda_\epsilon+\overline{\lambda}_\epsilon}\Psi_{\epsilon
1}(0)\mathrm{B}_\epsilon\varphi_{\epsilon1}(0)(
e^{-\lambda_\epsilon
r}-e^{-(2\lambda_\epsilon+\overline{\lambda}_\epsilon)r})-g_{\epsilon2,1}\Psi_{\epsilon
1}(0)\varphi_{\epsilon1}(0)+g_{\epsilon2,1},
\]
\[E_{\epsilon 1b}= \frac{-\overline{g}_{\epsilon1,2}}{2\lambda_\epsilon}\Psi_{\epsilon
1}(0)\mathrm{B}_\epsilon\overline{\varphi}_{\epsilon1}(0)(e^{-\overline{\lambda}_\epsilon
r}-e^{-(2\lambda_\epsilon+\overline{\lambda}_\epsilon)r})-\overline{g}_{\epsilon1,2}\Psi_{\epsilon
1}(0)\overline{\varphi}_{\epsilon1}(0),
\]
\[E_{\epsilon2}=-2 g_{\epsilon1,1}e^{-(2\lambda_\epsilon+\overline{\lambda}_\epsilon)r}
\Psi_{\epsilon 1}(0)\mathrm{B}_\epsilon\int_{-r}^0w_{\epsilon2,0}(\theta)e^{-(2\lambda_\epsilon+
\overline{\lambda_\epsilon})
\theta}d\theta-2g_{\epsilon1,1}\Psi_{\epsilon 1}(0)w_{\epsilon 2,0}(0),
\]
\[E_{\epsilon3}=-(g_{\epsilon2,0}+2\overline{g}_{\epsilon1,1})e^{-(2\lambda_\epsilon+\overline{\lambda}_\epsilon)r}
\Psi_{\epsilon 1}(0)\mathrm{B}_\epsilon\int_{-r}^0
w_{\epsilon1,1}(\theta)e^{-(2\lambda_\epsilon+\overline{\lambda}_\epsilon)
\theta}d\theta-\]
\[-(g_{\epsilon2,0}+2\overline{g}_{\epsilon1,1})\Psi_{\epsilon 1}(0)w_{\epsilon
1,1}(0),
\]
\[E_{\epsilon4}=-\overline{g}_{\epsilon0,2}e^{-(2\lambda_\epsilon+\overline{\lambda}_\epsilon)r}\Psi_{\epsilon
1}(0)\mathrm{B}_\epsilon\int_{-r}^0
w_{\epsilon0,2}(\theta)
e^{-(2\lambda_\varepsilon+\overline{\lambda}_\epsilon)
\theta}d\theta-\overline{g}_{\epsilon0,2}\Psi_{\epsilon 1}(0)w_{\epsilon0,2}(0).
\]

$\mathbf{E_{\epsilon1a}}.$
We have
\[E_{\epsilon1a}=-g_{\epsilon 21}\Psi_{\epsilon
1}(0)\mathrm{B}_\epsilon\varphi_{\epsilon1}(0)\int_{-r}^0e^{-\lambda_\epsilon(s+
r)}e^{(2\lambda_\epsilon+\overline{\lambda}_\epsilon)s}ds
-g_{\epsilon 21}\Psi_{\epsilon 1}(0)\varphi_{\epsilon1}(0) +g_{\epsilon 21}.
\]
We consider the function $\rho_{\epsilon}(s)=\varphi_{\epsilon
1}(0)e^{(2\lambda_\epsilon+\overline{\lambda}_\epsilon)s},\,s\in[-r,0]$
and remark that
\[E_{\epsilon1a}=g_{\epsilon 21}(1-\langle \Psi_{\epsilon1}\,,\,\rho_{\epsilon}  \rangle).
\]

But we know that $1=\langle  \Psi_{\epsilon1}, \varphi_{\epsilon
1}\rangle,$ and, substituting this above, we find
\[E_{\epsilon1a}=g_{\epsilon21}(\langle  \Psi_{\epsilon1}, \varphi_{\epsilon1}\rangle-\langle
\Psi_{\epsilon1}, \rho_\epsilon \rangle)=
g_{\epsilon21}\langle  \Psi_{\epsilon1}, \varphi_{\epsilon1}-\rho_\epsilon\rangle,
\]
and
\[\varphi_{\epsilon1}(s)-\rho_\epsilon(s)=\varphi_{\epsilon1}(0)\left(e^{\lambda_\epsilon
s}-e^{(2\lambda_\epsilon+\overline{\lambda}_\epsilon)s}\right)=e^{\lambda_\epsilon s}\left(1-
e^{2\mu_\epsilon s}\right)\varphi_{\epsilon1}(0)= \]
\begin{equation}\label{imp-rel}=-2\mu_\epsilon s\, e^{\lambda_\epsilon
s}\left(1+\frac{2\mu_\epsilon s}{2!}+\frac{(2\mu_\epsilon
s)^2}{3!}+...\right)\varphi_{\epsilon1}(0).
\end{equation}
The series in the paranthesis converges uniformly and its sum is a bounded
function on $[-r,0].$  We see that
\[\lim_{\epsilon\rightarrow 0}E_{\epsilon1a}/\mu_\epsilon=g_{21}\langle \Psi_{1}, \,\rho
\rangle,
\]
where $\rho(s)=-2se^{\lambda s}\varphi_{1}(0).$\\

Now we pass to $\mathbf{E_{\epsilon1b}}$.
\[E_{\epsilon1b}=
\frac{-\overline{g}_{\epsilon1,2}}{2\lambda_\epsilon}\Psi_{\epsilon1}(0)\mathrm{B}_\epsilon\overline{\varphi}_{\epsilon1}(0)(e^{-\overline{\lambda}_\epsilon
r}-e^{-(2\lambda_\epsilon+\overline{\lambda}_\epsilon)r})-
\overline{g}_{\epsilon1,2}\Psi_{\epsilon1}(0)\overline{\varphi}_{\epsilon1}(0)=
\]
\[=-
\overline{g}_{\epsilon1,2}\left[\Psi_{\epsilon1}(0)\overline{\varphi}_{\varepsilon1}(0)+\Psi_{\epsilon1}(0)\mathrm{B}_\epsilon\overline{\varphi}_{\varepsilon1}(0)e^{-(2\lambda_\epsilon+\overline{\lambda}_\epsilon)r}\frac{1}{-2\lambda_\epsilon}(1-e^{2\lambda_\epsilon 
r})\right]=
\]
\[=-
\overline{g}_{\epsilon1,2}\left[\Psi_{\epsilon1}(0)\overline{\varphi}_{\varepsilon1}(0)+\Psi_{\epsilon1}(0)\mathrm{B}_\epsilon\overline{\varphi}_{\varepsilon1}(0)\int_{-r}^0e^{-(2\lambda_\epsilon+\overline{\lambda}_\epsilon)(s+r)}e^{\overline{\lambda}_\epsilon 
s}ds\right].
\]
We define the function
$\eta_{\epsilon}(\zeta)=\Psi_{\epsilon1}(0)e^{-(2\lambda_\epsilon+\overline{\lambda}_\epsilon)\zeta},\,\zeta\in
[0,r]$ and get
\[E_{\epsilon1b} =-\overline{g}_{\epsilon1,2}\langle \eta_{\epsilon},\,\varphi_{\epsilon
2}\rangle.
\]

Since $\langle \Psi_{\epsilon1},\,\varphi_{\epsilon 2} \rangle=0,$ we have
$$-\langle \eta_{\epsilon},\,\varphi_{\epsilon 2}\rangle=\langle
\Psi_{\epsilon1},\,\varphi_{\epsilon 2} \rangle-\langle \eta_{\epsilon},\,\varphi_{\epsilon 2}
\rangle=\langle\Psi_{\epsilon1}-\eta_{\epsilon},\,\varphi_{\epsilon 2} \rangle,$$ and
\[\Psi_{\epsilon1}(\zeta)-\eta_{\epsilon}(\zeta)=\Psi_{\epsilon1}(0)\left(e^{-\lambda_\epsilon
\zeta}-e^{-(2\lambda_\epsilon+\overline{\lambda}_\epsilon)(\zeta)}
\right)=\Psi_{\epsilon1}(0)e^{-\lambda_\epsilon
\zeta}(1-e^{-(\lambda_\epsilon+\overline{\lambda}_\epsilon)\zeta})=
\]
\[=\Psi_{\epsilon1}(0)e^{-\lambda_\epsilon \zeta}2\mu_\epsilon \zeta(1-\frac{2\mu_\epsilon
\zeta}{2!}+\frac{(2\mu_\epsilon \zeta)^2}{3!}+...).
\]
Hence
\[\lim_{\epsilon\rightarrow 0}E_{\epsilon1b}/\mu_\epsilon= \overline{g}_{\epsilon1,2}
\langle \eta,\varphi_{\epsilon 2} \rangle,
\]
where
$\eta(\zeta)=2\zeta e^{-\omega i\zeta}\Psi_{1}(0),\,\,\zeta\in[0,r].$

\noindent$\mathbf{E_2}\,-\,\mathbf{E_4}.$ We
define
\[\alpha_{2,0}:=-2g_{\epsilon1,1},\,\,\alpha_{1,1}:=-g_{\epsilon2,0}-2\overline{g}_{\epsilon1,1},\,
\,\alpha_{0,2}:=-\overline{g}_{\epsilon0,2}.\]

The expressions $E_i,\, i=2,3,4, $ can be written as
\[\alpha_{j,k}\left(\Psi_{\epsilon 1}(0)\cdot w_{\epsilon j,k}(0)+ \Psi_{\epsilon 1}(0)\cdot
\mathrm{B}_\epsilon \int_{-r}^0 e^{-(2\lambda_\epsilon+\overline{\lambda}_\epsilon)(s+r)}
w_{\epsilon j,k}(s)ds     \right)=\alpha_{j,k}\langle \eta_\epsilon, w_{\epsilon j,k}\rangle,\;\;
\]
where $j,k>0,\,j+k=2.$

We know that
\[\langle \Psi_{\epsilon1}, w_{\epsilon j,k} \rangle =0,
\]
and we may write
\[\langle \eta_\epsilon, w_{\epsilon j,k}\rangle=
\langle \eta_\epsilon, w_{\epsilon j,k}\rangle-\langle
\Psi_{\epsilon1}, w_{\epsilon j,k} \rangle= \langle
\eta_\epsilon-\Psi_{\epsilon1}, w_{\epsilon j,k}\rangle,
\]
and
\[\eta_\epsilon(\zeta)-\Psi_{\epsilon1}(\zeta)=\Psi_{\epsilon
1}(0)(e^{-(2\lambda_\epsilon+\overline{\lambda}_\epsilon)\zeta}-e^{-\lambda_\epsilon\zeta})=\Psi_{\epsilon
1}(0)e^{-\lambda_\epsilon\zeta}\left[ e^{-2\mu_\epsilon\zeta}-1 \right]=
\]
\[=-2\mu_\epsilon\zeta e^{-\lambda_\epsilon\zeta}\left(
1-\frac{2\mu_\epsilon\zeta}{2!}+\frac{(2\mu_\epsilon\zeta)^2}{3!}-... \right)\Psi_{\epsilon
1}(0).
\]
Hence
\[\lim_{\epsilon \rightarrow 0} \langle \eta_\epsilon,\,w_{\epsilon j,k} \rangle/\mu_\epsilon =
-\langle \eta,\,w_{j,k} \rangle.
\]

\vspace{0.3cm}

From the above equalities we see that
\[\lim_{\epsilon\rightarrow
0}(E_{\epsilon2}+E_{\epsilon3}+E_{\epsilon4})/\mu_\epsilon=2g_{1,1}\langle\eta ,w_{2,0}\rangle+
(g_{2,0}+2\overline{g}_{1,1})\langle\eta ,w_{1,1}\rangle+\overline{g}_{0,2}\langle\eta
,w_{0,2}\rangle.
\]

\vspace{0.2cm}

Finally, we have, by setting $h_2(\epsilon)=\Psi_{\epsilon1}(0)\cdot
(\mathrm{B}_{\epsilon}R_{\epsilon1}-R_{\epsilon2})/\mu_\epsilon, $
$$\lim_{\epsilon \rightarrow 0}h_2(\epsilon) =
g_{21}\langle \Psi_1, \rho \rangle+\overline{g}_{1,2}\langle\eta
,\varphi_2\rangle+2g_{1,1}\langle\eta ,w_{2,0}\rangle+ (g_{2,0}+2\overline{g}_{1,1})\langle\eta
,w_{1,1}\rangle+\overline{g}_{0,2}\langle\eta ,w_{0,2}\rangle,
$$ where
\begin{eqnarray}
    \rho(s) &=& -2se^{\omega i s}\varphi_1(0),\,\,s\in[-r,0], \\
    \eta(\zeta) &=& 2\zeta e^{-\omega i\zeta}\Psi_1(0),\,\,\zeta\in[0,r].
\end{eqnarray}
This proves the second affirmation of the Proposition.$\Box$

\vspace{0.5cm}

The proof of Proposition 4.2 allows us to assert

\vspace{0.5cm}

\textbf{Theorem 4.1} \textit{The value in zero of the function $w_{2,1}$ is given by}
\[w_{2,1}(0)=x_1v_1+x_2v_2+...+x_nv_n,
\]
\textit{where} $v_1=\Psi_1(0)^T/|\Psi_1(0)|,$ \textit{the set of vectors} $v_1,...,v_n$ \textit{is
an orthonormal basis of} $\mathbb{C}^n$ \textit{and} $x_1,...,x_n$ \textit{are solutions of the
system}
$$\left\{\begin{array}{lll}
(\Psi_1(0)\widetilde{\mathrm{M}}v_1)x_1+(\Psi_1(0) \widetilde{\mathrm{M}}v_2)x_2+...+(\Psi_1(0)
\widetilde{\mathrm{M}}v_n)x_n&=& \widetilde{R},\\
(v_2\cdot \mathrm{M}v_1)x_1+(v_2\cdot \mathrm{M}v_2)x_2+...+(v_2\cdot
\mathrm{M}v_n)x_n&=&v_2\cdot(\mathrm{B}R_1-R_2),\\
                                                      &\vdots &                     \\
(v_n\cdot \mathrm{M}v_1)x_1+(v_n\cdot \mathrm{M}v_2)x_2+...+(v_n\cdot
\mathrm{M})v_nx_n&=&v_n\cdot(\mathrm{B}R_1-R_2),\\
\end{array}\right. $$
\textit{with} $\widetilde{M}=2 \left(\mathrm{B} e^{-\omega i r} r+\mathrm{I}\right),$
$\,\,\,M=\omega i\mathrm{I}-\mathrm{A}-\mathrm{B}e^{-\omega ir},$
\[\widetilde{R}=g_{21}\langle \Psi_1, \rho \rangle+\overline{g}_{1,2}\langle\eta
,\varphi_2\rangle+2g_{1,1}\langle\eta ,w_{2,0}\rangle+(g_{2,0}+2\overline{g}_{1,1})\langle\eta
,w_{1,1}\rangle+\overline{g}_{0,2}\langle\eta ,w_{0,2}\rangle,
\]
$\Psi_1(\zeta)=[\psi_1(0)\varphi_1(0)+\psi_1(0)
 \mathrm{B}\varphi_1(0)e^{-\omega ri}r]^{-1}\,e^{-\omega i \zeta}\psi_1(0)$, \\

\noindent $ \rho(s) = -2se^{\omega i s}\varphi_1(0),\,\,s\in[-r,0],$ \textit{and} $\eta(\zeta) =
2\zeta e^{-\omega i\zeta}\Psi_1(0),\,\,\zeta\in[0,r].$\\

\vspace{0.3cm}

Once $w_{2,1}(0)$ known, $w_{2,1}(-r)$ is computed from \eqref{alg-sys1}.

\section{Appendix}

\textbf{I.} We compute here the matrix $E$ with elements $e_{ij}=\langle \psi_i, \,\varphi_j
\rangle$.

\[e_{11}=\langle \psi_1, \,\varphi_1 \rangle=\psi_1(0)\varphi_1(0)+
\psi_1(0)\mathrm{B}\varphi_1(0)\int_{-r}^0e^{-(\zeta+r)\omega i}e^{\zeta
\omega i}d\zeta=\]
\[=\psi_1(0)\varphi_1(0)+\psi_1(0)\mathrm{B}\varphi_1(0)e^{-\omega
ri}r.
\]

\[e_{21}=\langle \psi_2, \,\varphi_1 \rangle=\psi_2(0) \varphi_1(0)+
\psi_2(0) \mathrm{B}\varphi_1(0)\int_{-r}^0e^{(\zeta+r)\omega i}e^{\zeta
\omega
i}d\zeta=\]
\[=\psi_2(0)\varphi_1(0)+\psi_2(0) \mathrm{B}\varphi_1(0)e^{\omega
ri}\frac{1}{2\omega i}(1-e^{-2\omega ir})=
\]
\[=\psi_2(0)\varphi_1(0)+\frac{1}{2\omega i}(\psi_2(0)  \mathrm{B}e^{\omega
ri}\varphi_1(0)-\psi_2(0) \mathrm{B} e^{-\omega
ri}\varphi_1(0)).
\]
In the expression above, we use \eqref{ev-adj} to eliminate $\psi_2(0) Be^{\omega
ri}\varphi_1(0)$ and \eqref{ev} to eliminate $\psi_2(0) \mathrm{B} e^{-\omega
ri}\varphi_1(0)$. We get
\[e_{21}=\psi_2(0)\varphi_1(0)+\frac{1}{2\omega i}\left(-\psi_2(0)\omega
i\varphi_1(0)-\psi_2(0)\mathrm{A}\varphi_1(0)-\psi_2(0)\omega
i\varphi_1(0)+\psi_2(0)A\varphi_1(0)\right)=\]
\[=\psi_2(0)\varphi_1(0)-\psi_2(0)\varphi_1(0)=0.
\]
For $e_{12}$, as for $e_{21},$ we find
\[e_{12}=0.
\]
Then
\[e_{22}=\langle \psi_2, \,\varphi_2 \rangle=\psi_2(0)\varphi_2(0)+
\psi_2(0)\mathrm{B}\varphi_2(0)\int_{-r}^0e^{(\zeta+r)\omega i}e^{-\zeta
\omega i}d\zeta=\]
\[=\psi_2(0)\varphi_2(0)+\psi_2(0)\mathrm{B}\varphi_2(0)e^{\omega
ri}r,
\]

Now, by using
\[\left(%
\begin{array}{c}
  \Psi_1 \\
  \Psi_2 \\
\end{array}%
\right)=E^{-1}\left(%
\begin{array}{c}
  \psi_1 \\
  \psi_2 \\
\end{array}%
\right)=\frac{1}{\det E}\left(
                         \begin{array}{cc}
                            e_{22} & 0\\
                           0  &  e_{11} \\
                         \end{array}
                       \right)\left(%
\begin{array}{c}
  \psi_1 \\
  \psi_2 \\
\end{array}\right)
\]
we will get $\Psi_1(\zeta)=\Psi_1(0)e^{-\omega i \zeta}$, with
$$\Psi_1(0)=\frac{e_{22}}{\det\,
E}\psi_1(0)=\frac{1}{e_{11}}\psi_1(0)=\frac{1}{\psi_1(0)\varphi_1(0)+\psi_1(0)
 \mathrm{B}\varphi_1(0)e^{-\omega
ri}r}\psi_1(0).$$

\vspace{0.5cm}

\textbf{II.} We compute the matrix $E_\epsilon$ with elements $\widetilde{e}_{ij}=\langle
\psi_{\epsilon i}, \,\varphi_{\epsilon j} \rangle$:
$$\widetilde{e}_{11}=\langle \psi_{\epsilon1},\varphi_{\epsilon1}
\rangle=\psi_{\epsilon1}(0)\varphi_{\epsilon1}(0)+\psi_{\epsilon1}(0)\mathrm{B}_\epsilon\varphi_{\epsilon1}(0)\int_{-r}^0e^{-\lambda_\epsilon(s+r)}e^{\lambda_\epsilon 
s}ds=
$$
$$=\psi_{\epsilon1}(0)\varphi_{\epsilon1}(0)+\psi_{\epsilon1}(0)
\mathrm{B}_\epsilon\varphi_{\epsilon1}(0)e^{-\lambda_\epsilon r}r.
$$
$$
\widetilde{e}_{12}=\langle \psi_{\epsilon1},\varphi_{\epsilon2}
\rangle=\psi_{\epsilon1}(0)\varphi_{\epsilon2}(0)+\psi_{\epsilon1}(0)\mathrm{B}_\epsilon\varphi_{\epsilon2}(0)\int_{-r}^0e^{-\lambda_\epsilon(s+r)}e^{\overline{\lambda}_\epsilon 
s}ds=
$$
\[=\psi_{\epsilon1}(0)\varphi_{\epsilon2}(0)+\frac{1}{-\lambda_\epsilon
+\overline{\lambda_\epsilon}}\psi_{\epsilon1}(0)(\mathrm{B}_\epsilon\varphi_{\epsilon2}(0)e^{-\lambda_\epsilon
r}-\mathrm{B}_\epsilon\varphi_{\epsilon2}(0)e^{-\overline{\lambda_\epsilon} r})=
\]
\[=\psi_{\epsilon1}(0)\varphi_{\epsilon2}(0)+\frac{1}{-\lambda_\epsilon
+\overline{\lambda}_\epsilon}\left[(\lambda_\epsilon\psi_{\epsilon1}(0)
-\psi_{\epsilon1}(0)\mathrm{A}_\epsilon)\varphi_{\epsilon2}(0)-
\psi_{\epsilon1}(0)(\overline{\lambda_\epsilon}\varphi_{\epsilon2}(0)-\mathrm{A}_\epsilon\varphi_{\epsilon2}(0))\right]=
\]
\[=\psi_{\epsilon1}(0)\varphi_{\epsilon\,2}(0)+\frac{1}{-\lambda_\epsilon
+\overline{\lambda}_\epsilon}(\lambda_\epsilon-\overline{\lambda}_\epsilon)\psi_{\epsilon1}(0)\varphi_{\epsilon2}(0)=0.
\]
Here we used the relations
\[\psi_{\epsilon1}(0)(-\lambda_\epsilon \mathrm{I}+\mathrm{A}_\epsilon+\mathrm{B}_\epsilon
e^{-\lambda r})=0, \,\, (\overline{\lambda}_\epsilon
\mathrm{I}-\mathrm{A}_\epsilon-\mathrm{B}_\epsilon
e^{-\overline{\lambda_\epsilon}r})\varphi_{\epsilon2}(0)=0.
\]
Now, by taking the complex conjugate we obtain
\[\widetilde{e}_{21}=\langle \psi_{\epsilon2}(0),\,\varphi_{\epsilon1}(0) \rangle =0,
\]
\[\widetilde{e}_{22}=\langle \psi_{\epsilon2},\varphi_{\epsilon2}
\rangle=\psi_{\epsilon2}(0)\cdot\varphi_{\epsilon2}(0)+\psi_{\epsilon2}(0)\mathrm{B}_{\epsilon}\varphi_{\epsilon2}(0)\int_{-r}^0e^{-\overline{\lambda}_\epsilon(s+r)}e^{\overline{\lambda}_\epsilon 
s}ds=\]
\[=\psi_{\epsilon2}(0)\cdot\varphi_{\epsilon2}(0)+\psi_{\epsilon2}(0)\mathrm{B}_\epsilon\varphi_{\epsilon2}(0)e^{-\overline{\lambda}_\epsilon 
r}r.
\]

We have $\det \,E=\widetilde{e}_{11}\widetilde{e}_{22}$ and
\[\Psi_{\epsilon
1}(\zeta)=\frac{\widetilde{e}_{22}}{\widetilde{e}_{11}\widetilde{e}_{22}}\psi_{\epsilon
1}(\zeta)=\frac{1}{\widetilde{e}_{11}}e^{-\lambda_\epsilon \zeta}\psi_{\epsilon 1}(0).
\]

\end{document}